\documentclass[11pt,leqno]{article} 
\usepackage{graphics}
\newtheorem{thm}{Theorem}[section]
\newtheorem{lma}{Lemma}[section]

\newcommand{\beqa}{\begin{eqnarray}}
\newcommand{\eeqa}{\end{eqnarray}}

\newcommand{\pf}{\noindent {\bf Proof:} $\s$ }
\newcommand{\epf}{ \hfill$\diamondsuit$ \medskip}
\newcommand{\md}{\medskip}

\newcommand{\beq}{\begin{equation}}
\newcommand{\eeq}{\end{equation}}
\newcommand{\lbl}{\label}
\newcommand{\s}{\; \;}
\newcommand{\noi}{\noindent}

\newcommand{\p}{\varphi}
\newcommand{\la}{\lambda}
\newcommand{\mb}{\mbox}
\newcommand{\ra}{\rightarrow}
\newcommand{\al}{\alpha}

\title{Existence and uniqueness of solutions for a class of $p$-Laplace equations on a ball}

\author{
Philip Korman   \\ 
Department of Mathematical Sciences \\ 
University of Cincinnati \\ 
Cincinnati Ohio 45221-0025 \\
}

\date{}

\begin{document}

\maketitle
\begin{abstract} 
\noindent
For a class of equations generalizing the model case 
\[
\Delta _p u-a(r)u^{p-1}+b(r)u^q=0 \s \mbox{in $B$}, \s u=0 \s \mbox{on $\partial B$},
\]
where $B$ is the unit ball in $R^n$, $n \geq 1$, $r=|x|$, $p,q>1$,  and $\Delta _p$ denotes the $p$-Laplace operator, we give conditions for the existence and uniqueness of positive solution. In case $n=1$, we give a  more general result.
 \end{abstract}

\begin{flushleft}
\end{flushleft}

\begin{flushleft}
\end{flushleft}

\section{Introduction}
\setcounter{equation}{0}
\setcounter{thm}{0}
\setcounter{lma}{0}

We study the existence and  uniqueness of positive solutions of the Dirichlet problem
\beq
\lbl{i1}
\Delta _p u+f(r,u)=0 \s \mbox{in $B$}, \s u=0 \s \mbox{on $\partial B$},
\eeq
where $B$ is the unit ball in $R^n$, $n \geq 1$, $r=|x|$, and $\Delta _p$ denotes the $p$-Laplace operator
\[
\Delta _p u=\mbox{div} \left( |\nabla u|^{p-2} \nabla u \right), \s p>1 \,.
\]
We shall assume that $f_r(r,u) \leq 0$ for all $r \in [0,1]$, and $u>0$, which implies (in case $1<p \leq 2$) that positive  solutions are radially symmetric, i.e., $u=u(r)$, 
and the equation (\ref{i1}) becomes
\[
\left(u'|u'|^{p-2} \right)'+\frac{n-1}{r} u'|u'|^{p-2} +f(r,u)=0, \s 0<r<1, \s u'(0)=u(1)=0 \,.
\]
In case $p=2$, exact multiplicity results were given in P. Korman, Y. Li and T. Ouyang \cite{KLO}, T. Ouyang and J. Shi \cite{OS}, and P. Korman \cite{K}, and other papers of the same authors. Uniqueness results, in case $p=2$, were studied by many authors, beginning with B. Gidas, W.-M. Ni and L. Nirenberg \cite{GNN}, for $f(u)=u^q$, see the references in \cite{KLO}, \cite{OS}. When $p \ne 0$, there are a number of technical difficulties, some of which were overcome only recently, through the efforts of a number of people. In particular, F. Brock \cite{B} and L. Damascelli and F.  Pacella \cite{DP} have extended the symmetry results of B. Gidas, W.-M. Ni and L. Nirenberg \cite{GNN}. J. Serrin and H. Zou \cite{SZ} have proved a 
Liouville type result, which  implies a priori estimates of B. Gidas and J. Spruck type \cite{GS}, provided that the point of maximum of solution is bounded away from the boundary, as was observed in A. Aftalion and F. Pacella \cite{AP1}, and J. Fleckinger and W. Reichel \cite{FR}. A. Aftalion and F. Pacella \cite{AP1} pointed out the appropriate function spaces, in order to apply the implicit function theorem. In fact,  A. Aftalion and F. Pacella \cite{AP1} have proved uniqueness of solution to (\ref{i1}), for $f(r,u)$ similar to ours.  In this paper, we provide a simpler proof of the crucial step in \cite{AP1}, involving the non-degeneracy of solutions. Moreover, we clarify the proper conditions, and provide the existence part for the model examples.
\md 

Recall that solution of (\ref{i1}) is called non-degenerate if the corresponding linearized problem (the problem (\ref{31}) below) has only the trivial solution. Non-degeneracy results were usually proved by  using the method of test functions, see P. Korman, Y. Li and T. Ouyang \cite{KLO}, T. Ouyang and J. Shi \cite{OS}, and the references in those papers. A variation on that approach, involving maximum principle, was used in A. Aftalion and F. Pacella \cite{AP1}. Here we use a simpler method, based on some identities of M. Tang, and his method \cite{T}. We have already used this approach, see \cite{K1}, in case $p=2$ (M. Tang's  results are for $p=2$ case). In this paper, we extend M. Tang's \cite{T} identities  for $p \ne 2$ case, allowing a considerable streamlining of the proof of non-degeneracy. These new identities are likely to be useful for other problems.
\md 

Our main example is $f(r,u)=-a(r)u^{p-1}+b(r)
u^q$, with $p,q >1$, and  $p-1<q$. This function is negative for small $u$, and positive for $u$ large. We give conditions for existence and uniqueness of positive solution in this case. We show that in one space dimension, and for $f=f(u)$, the same result remains true, if $f(u)$ is arbitrarily modified in the region where it is negative. This surprising result was proved, in case $p=2$, by R. Schaaf \cite{S}. The uniqueness part was proved (for any  $p >1$) by J. Cheng \cite{C}. Both of these works were using time maps, while we are using a more flexible bifurcation approach.

\section{Existence and uniqueness for a class of non-autonomous problems}
\setcounter{equation}{0}
\setcounter{thm}{0}
\setcounter{lma}{0}
We consider positive solutions of the problem
\beq
\lbl{20}
\Delta _p u+f(|x|,u)=0, \s \mb{for $|x|<1$}, \s u=0 \s \mb{on $|x|=1$}.
\eeq
We assume that $f(r,u) \in C^1([0,1] \times \bar R_+)$ is a continuously differentiable function (with $r=|x|$), and 
\beq
\lbl{21}
f(r,0)=0, \s \mb{and} \s f_r(r,u) \leq 0, \s \s  \mb{for all $r \in [0,1]$, and $u>0$} \,.
\eeq
If $1< p \leq 2$, then  in view of  B. Gidas, W.-M. Ni and L. Nirenberg \cite{GNN} and L. Damascelli and F.  Pacella \cite{DP}, positive solutions of (\ref{20}) are radially symmetric,  and hence they satisfy
\beq
\lbl{30}
\p(u'(r))'+\frac{n-1}{r} \p(u'(r))+f(r,u(r))=0, \s 0<r<1, \s u'(0)=u(1)=0 \,.
\eeq
where we denote $\p (t)=t|t|^{p-2}$. In case $p>2$, we restrict our attention to the radial solutions of (\ref{20}), i.e., we again consider (\ref{30}). Following B. Franchi et al \cite{F}, we consider {\em classical} solutions of (\ref{30}), i.e., we assume that $u \in C^1[0,1]$ and $\p(u') \in C^1(0,1]$. As a consequence,  if $u'(r_1)=0$ for some $r_1 \in (0,1)$, we may define by continuity $\p (u'(r_1))=0$. 
\md 

The following lemma gives a known condition for the Hopf's lemma to hold, see J.L. V\'{a}zquez \cite{V}, or A. Aftalion and F. Pacella \cite{AP1}. 

\begin{lma}\lbl{lma:00}
Assume that $f(r,u) \in C^1([0,1] \times \bar R_+)$ satisfies $f(r,0)=0$, for all $r \in [0,1]$, and either $f$ is non-negative, or for some $c_0>0$
\beq
\lbl{14-}
-f(r,s)<c_0 s^{p-1}, \s \mbox{for small $s>0$, and $r$ near $1$} \,.
\eeq
If $u(r)$ is a positive solution of (\ref{30}), then
\beq
\lbl{14a}
u'(1)<0 \,.
\eeq
\end{lma}

The next lemma extends the Proposition 1.2.6 in B. Franchi et al \cite{F}, which considered the case of $f=f(u)$, i.e., independent
of $r$. Define $F(r,u)=\int _0^u f(r,t) \, dt$.

\begin{lma}\lbl{lma:01}
Assume that
\beq
\lbl{14b} 
F_r(r,u) \leq 0, \s \s  \mb{for all $r \in [0,1]$, and $u>0$} \,;
\eeq
\beqa
\lbl{14c}
& \mbox{for each $r_0 \in [0,1)$, the function $f(r_0,u)$ is either positive for all $u>0$,} \\
& \mbox{or it changes sign once from negative, for small $u$, to positive. }       \nonumber
\eeqa
Then any positive solution of (\ref{30}) satisfies
\beq
\lbl{14d}
u'(r)<0 \s\s \mbox{ for all $r \in (0,1)$} \,,
\eeq
\beq
\lbl{14e}
f(0,u(0))>0 \,.
\eeq
\end{lma}

\pf
We show first that $u(r)$ has no local minimums. Indeed, let $r_0 \in [0,1]$ be a point of local minimum, i.e., $u'(r_0)=0$, $u''(r_0) \geq 0$. We claim that 
\beq
\lbl{14f}
f(r_0,u(r_0)) \leq 0 \,.
\eeq
If $p \geq 2$, this follows by evaluating the equation (\ref{30}) at $r_0$. In general, we can find a sequence $\{ r_n \}$ tending to $r_0$ from the right, at which $u'(r_n) > 0$, $u''(r_n) \geq 0$. Evaluating the equation (\ref{30}) at $r_n$, we have $f(r_n,u(r_n)) < 0$, and taking the limit we conclude (\ref{14f}). In case $f$ is positive, (\ref{14f}) implies a contradiction. In the other case, (\ref{14f}) implies that $F(r_0,u(r_0))<0$, and hence the ``energy" $E(r) \equiv \frac{p-1}{p} |u'(r)|^p + F(r,u(r))$ satisfies $E(r_0)<0$.
We have 
\[
E'(r)=-\frac{n-1}{r} \p(u')u'+F_r(r,u) \leq 0 \,,
\]
and so the energy is non-increasing. But $E(1)=\frac{p-1}{p} |u'(1)|^p \geq 0$, a contradiction. We conclude that $u'(r) \leq 0$ for all $r$. To get the strict inequality, assume that $u'(r_1)=0$ at some $r_1 \in (0,1)$. As discussed above, we have $\p (u'(r_1))=0$, and we also have $\frac{d}{dr} \p (u'(r_1))=0$, since $r_1$ is not an extremum. From  the equation (\ref{30}), $f(r_1,u(r_1)) = 0 $, implying that $E(r_1)<0$, which results in the same contradiction as before. Finally, we conclude (\ref{14e}), by using the same argument one more time, since the opposite inequality would imply
$E(0)<0$.
\epf

We assume that $f(r,u)$ satisfies the following additional conditions, which  are similar to the ones in \cite{AP1}
\beq
\lbl{22}
uf_u(r,u)-(p-1)f(r,u)>0 \s \s  \mb{for all $r \in [0,1]$, and $u>0$} \,;
\eeq
For any positive solution of (\ref{30})
\beqa
\lbl{23}
& \al (r)=\frac{pf(r,u(r))+rf_r(r,u(r))}{u(r)f_u(r,u(r))-(p-1)f(r,u(r))} \\
& \mb{is a non-increasing function of $r$, for $r \in (0,1)$}. \nonumber
\eeqa

\md

We shall need to consider the  linearized problem corresponding to (\ref{30}) (here $w=w(r)$)
\beqa
\lbl{31}
& \left(\p '(u'(r))w'(r) \right)'+\frac{n-1}{r} \p'(u'(r))w'(r)+f_u(r,u(r))w(r)=0, \s 0<r<1, \\
& w'(0)=w(1)=0 \,. \nonumber
\eeqa
We will show that under the above  conditions any positive solution of (\ref{30}) is non-degenerate, i.e., the problem (\ref{31}) admits only the trivial solution.  The following  technical lemma we proved in \cite{K1}. We include its proof for completeness.
\begin{lma}\lbl{lma:20}
Let $u(r)$ be a positive solution of (\ref{30}), and assume that the function $f(r,u)$ satisfies the conditions (\ref{21}) and (\ref{22}). Then the function $f(r,u(r))$ can change sign at most once on $(0,1)$.
\end{lma}

\pf
Let $\xi \in (0,1)$ be such that $f(\xi,u(\xi))=0$. We claim that $f(r,u(r))>0$ for all $r \in [0,\xi)$. Indeed, by (\ref{22}) we conclude that $f_u(\xi,u(\xi))>0$, and in general $f_u(r,u(r))>0$, so long as $f(r,u(r))>0$, and hence
\[
\frac{d}{dr} f(r,u(r))=f_r(r,u(r))+ f_u(r,u(r))u'(r)<0,
\]
and the claim follows. So that, if the function $f(r,u(r))$ is positive near $r=1$, it is positive for all $r \in [0,1)$. If, on the other hand, $f(r,u(r))$ is negative near $r=1$, it will change sign exactly once on $[0,1)$ (since it cannot stay negative for all $r$, by the maximum principle).
\epf

This lemma implies that either $f(r,u(r)) > 0$ for all $r \in [0,1)$, or else there is a $r_2 \in (0,1)$, so that $f(r,u(r))>0$ on $[0,r_2)$ and $f(r,u(r))<0$ on $(r_2,1)$. We put these cases together, we defining $r_2=1$ in the first case. I.e., the last lemma implies that
\beqa
\lbl{*}
& \mbox{there is a $r_2 \in (0,1]$, so that $f(r,u(r))>0$ on $[0,r_2)$} \\
& \mbox{ and, in case $r_2<1$, we have: $f(r,u(r))<0$ on $(r_2,1)$} \,. \nonumber
\eeqa
\md 

We shall consider the following two functions, depending on the solutions of (\ref{30}) and (\ref{31}):
\beq
\lbl{32}
\xi(r)=r^{n -1} \left[ (p-1) \p(u'(r))w(r)-\p'(u'(r))u(r)w'(r) \right] \,;
\eeq
\beq
\lbl{33}
\s\s\s\s T (r)=r^n \left[ (p-1) \p(u'(r))w'(r)+f(r,u(r)) w(r) \right]+(n-p)r^{n-1} \p(u'(r))w(r) \,.
\eeq
In case $p=2$, these functions were introduced by M. Tang \cite{T}. The following crucial lemma is proved by a direct computation.
\begin{lma}\lbl{lma:30}
For any solutions of  (\ref{30}) and (\ref{31}), we have
\beq
\lbl{35}
\xi'(r)=r^{n -1} w(r) \left[ u(r)f_u(r,u(r))-(p-1)f(r,u(r)) \right] \,;
\eeq
\beq
\lbl{36}
T '(r)=r^{n -1} w(r) \left[ p f(r,u(r))+r f_r (r,u(r)) \right] \,.
\eeq
\end{lma}

We shall need  the following functions,  depending on the solutions of (\ref{30}), see  \cite{OS}, \cite{T},  \cite{AP1}, and \cite{K1}.
\beq
\lbl{38}
\s\s\s\s Q(r)=r^n \left[ (p-1) \p \left( u'(r) \right)u'(r)+u(r)f(r,u(r)) \right] +(n-p)r^{n-1} \p \left( u'(r) \right)u(r);
\eeq
\beq
\lbl{39}
\s\s\s\s P(r)=r^n \left[ (p-1) \p \left( u'(r) \right)u'(r)+pF(r,u(r)) \right] +(n-p)r^{n-1} \p \left( u'(r) \right)u(r) \,,
\eeq
where we again denote $F(r,u)=\int _0^u f(r,t) \,dt$.
Observe that 
\beq
\lbl{39.1}
P(0)=0, \s \mbox{and} \s P(1)=\p \left( u'(1) \right) u'(1)>0 \,,
\eeq
provided (\ref{14-}) holds, and
\beqa 
\lbl{40}
&  P'(r) = r^{n-1} \left[npF(r,u(r))-(n-p) u(r)f(r,u(r))+pr F_r(r,u(r)) \right] \\
& \equiv  r^{n-1}I(r). \nonumber
\eeqa

We shall need the following lemma.
\begin{lma}\lbl{lma:31}
Let $u(r)$ be a positive solution of (\ref{30}), and assume that the function $f(r,u)$ satisfies the condition (\ref{21}). Then 
any solution of the linearized problem (\ref{31}), $w(r)$, cannot vanish in the region where $f(r,u(r))<0$ (i.e., on $(r_2,1)$, see (\ref{*})).
\end{lma}

\pf
Assume that the contrary is true, and let $\tau $ denote the largest root of $w(r)$ in the region where $f(r,u(r))<0$. We may assume that $w(r)>0$ on $(\tau,1)$. Then integrating the formula (\ref{36})  over $(\tau,1)$, we have
\[
(p-1) \p \left( u'(1) \right) w'(1)-\tau ^n (p-1) \p \left( u'(\tau) \right) w'(\tau)=\int_\tau^1 \left( pf(r,u(r))+rf_r(r,u(r)) \right) wr^{n-1} \, dr.
\]
We have a contradiction, since the left hand side is non-negative, while the integral on the right is negative. 
\epf

Next, we present the crucial non-degeneracy result.

\begin{thm}\lbl{thm:20}
Let  $u(r)$ be a positive solution of (\ref{30}). Assume that the conditions (\ref{21}),  (\ref{14-}),  (\ref{14b}),  (\ref{14c}), (\ref{22}) and  (\ref{23}) hold. In the case $p<n$, assume additionally the following two conditions:
\beq
\lbl{42}
u(r) f(r,u(r))-p F(r,u(r))>0 \s\s \mb{for $r \in (0,r_2)$} ;
\eeq
\beq
\lbl{42.1}
 \mbox{the function $I(r)=npF(r,u(r))-(n-p) u(r)f(r,u(r))+pr F_r(r,u(r))$} \,, 
\eeq
defined in (\ref{40}), satisfies any one of the following three conditions:

\noi
{\bf (i)}
 $I(r)>0$ on $(0,r_2)$,
\smallskip

\noi
{\bf (ii)} $I(r)<0$ on $(0,1)$,
\smallskip

\noi
{\bf (iii)} $I(r)>0$ on $(0,r_0)$ and $I(r)<0$ on $(r_0,1)$, for some $r_0 \in (0,1)$.
\smallskip

Then $u(r)$ is a non-degenerate solution, i.e., the  corresponding linearized problem (\ref{31}) admits only the trivial solution.
\end{thm}

\pf
With $r_2$ as defined by (\ref{*}), we claim that
\beq
\lbl{43}
Q(r)>0 \s \mbox{on $[0,r_2)$} \,.
\eeq
In the case $p \geq n$, this follows immediately from the definition of $Q(r)$ in (\ref{38}). In the other case, $p<n$,  observe that $P(r)>0$ on  $[0,r_2)$ by (\ref{39.1}) and (\ref{40}).
We  write
\[
Q(r)=P(r)+r^n \left[u(r) f(r,u(r))-p F(r,u(r)) \right] \,,
\]
 and hence $Q(r)>0$ on  $[0,r_2)$ by our condition (\ref{42}).
\md 

We now define the function $O(r)= \gamma \xi (r)- T(r)$, which in view of (\ref{35}) and (\ref{36}) satisfies
\beq
\lbl{44}
O'(r)= \left[u(r)f_u(r,u(r))-(p-1)f(r,u(r)) \right]w(r)r^{n-1} \left[ \gamma -\al (r) \right],
\eeq
with $\al (r)$ as defined by (\ref{23}), and $\gamma$ is a constant, to be selected. We may assume that $w(0)>0$. We claim that the function $w(r)$ cannot have any roots inside $(0,1)$.  Assuming otherwise, let $\tau _1$ be the smallest root of $w(r)$, i.e., $w(r)>0$ on $[0,\tau _1)$. Let $\tau _2 \in (\tau _1,1]$ denote the second root of $w(r)$. We now fix $\gamma=\al (\tau _1)$. By the monotonicity of $\al (r)$, the function $\gamma -\al (r)$ is non-positive on $[0,\tau _1)$ and non-negative  on $(\tau _1,\tau _2)$. Hence,  $O'(r) \leq 0$ (at $r=\tau _1$, both $w(r)$ and $\gamma -\al (r)$ change sign).  Since $O(0)=0$, we  conclude  that
\beq
\lbl{45}
O(r) \leq 0 \s \mb{for all $r \in [0,\tau _2]$} \,.
\eeq

There are two possibilities for the second root.
\smallskip

\noi
{\bf Case (i)} $\tau _2=1$, i.e., $w(r) < 0$ on $(\tau _1,1)$. From (\ref{45}) we have $O(1) \leq 0$. On the other hand, $O(1)=- T (1)=- (p-1) \p \left(u'(1) \right) w'(1)>0$, by (\ref{14a}), a contradiction.
\md

\noi
{\bf Case (ii)} $\tau _2<1$. Observe that $f(\tau _2,u(\tau _2))>0$. Indeed, assuming otherwise, we conclude by Lemma \ref{lma:20} that $f(r,u(r)) \leq 0$ on $(\tau _2,1)$. But this contradicts Lemma \ref{lma:31}.
 Applying Lemma \ref{lma:20} again, we conclude that $f(r,u(r))>0$ over $[0,\tau _2)$. It follows that $\tau _2 <r_2$, i.e., by (\ref{43}), $Q(r)>0$ on $[0,\tau _2)$.
\md

Since $\xi (\tau _1)>0$, while $\xi (\tau _2)<0$, we can find a point $t \in (\tau _1,\tau _2)$, such that $\xi (t)=0$, i.e.,
\beq
\lbl{46}
\frac{u(t)}{w(t)}=(p-1) \, \frac{\p \left(u'(t) \right)}{\p' \left(u'(t) \right)w'(t)}=\frac{u'(t)}{w'(t)} \,.
\eeq
Since $t \in (0,r_2)$, we have 
\beq
\lbl{48}
Q(t)>0.
\eeq
In view of (\ref{45}),
\[
T(t)=-O(t) \geq 0.
\]
On the other hand, using (\ref{46}) and (\ref{48}),
\[
T(t)=\left[ t^n \left((p-1) \p \left(  u' \right) w' \frac{u}{w}+f(t,u)u \right) +(n-p) t^{n-1} \p \left(  u' \right) u \right] \frac{w}{u}=Q(t)\frac{w(t)}{u(t)}<0,
\]
giving us a contradiction.
\md

It follows that $w(r)$ cannot have any roots, i.e., we may assume that $w(r)>0$ on $[0,1)$. But that is impossible, as can be seen by integrating (\ref{35}) over $(0,1)$. Hence $w \equiv 0$.
\epf

Our main example is   the problem (here $r=|x|$)
\beqa
\lbl{50}
& \p(u'(r))'+\frac{n-1}{r} \p(u'(r))  -a(r)u^{p-1}+b(r)
u^q=0, \s \s \mb{for} \s r \in (0,1), \\
& u =0 \s \mb{for $r=1$}. \nonumber
\eeqa
We consider the sub-critical case 
\beq
\lbl{49} 
\min (1,p-1)<q<\frac{np-n+p}{n-p} \,, 
\eeq
and assume that the functions $a(r)$, $b(r) \in C^1[0,1] $ satisfy
\beq
\lbl{51}
a(r)>0, \s b(r)>0, \s a'(r)>0, \s b'(r)<0 \s \s \mb{for } \, r \in (0,1) \,.
\eeq
 We define the functions
\[
A(r) \equiv pa(r)+ra'(r) \,,
\]
\[
B(r) \equiv \left(\frac{np}{q+1}-(n-p) \right)b(r)+\frac{pr b'(r)}{q+1} \,.
\]
Observe that $\frac{np}{q+1}-(n-p)>0$ for subcritical $q$, i.e., when (\ref{49}) holds.

\begin{thm}\lbl{thm:10}
In addition to the conditions (\ref{49}) and (\ref{51}), assume that the function $A(r)$ is positive and non-decreasing, while the function $B(r)$ is positive on $(0,1)$. Assume also that the functions $\frac{r b'(r)}{b(r)}$ and $r b'(r)$ are non-increasing on $(0,1)$.
Then any positive solution of the problem (\ref{50}) is non-degenerate.
\end{thm}

\pf
We shall verify the conditions of the Theorem \ref{thm:20}. We have
\[
uf-p F=b(r) u^{q+1} (1-\frac{p}{q+1})>0 \s\s \mb{for all $r \in [0,1)$} \,,
\]
verifying (\ref{42}).  Compute
\[
(q-(p-1)) \al (r)=- \frac{A(r)}{b(r)(u(r))^{q-p+1}}+p+  \frac{r b'(r)}{b(r)}.
\]
In view of our assumptions, $\al (r)$ is a non-increasing function, for all $r \in [0,1)$, verifying the condition (\ref{23}).
Compute
\[
I(r)=- A(r)u^p+B(r)u^{q+1}=B(r) u^p \left[-\frac{ A(r)}{ B(r)}+u^{q-p+1} \right].
\]
Since by our conditions, $B(r)$ is non-increasing, it follows that 
the quantity in the square bracket is a non-increasing function, which is negative near $r=1$. Hence, either $I(r)$ is negative over $(0,1)$, or it changes sign exactly once,  from positive to negative  thus verifying either part (ii), or part (iii), of the condition (\ref{42.1}).
Hence, the Theorem \ref{thm:20} applies, implying that any positive solution of the problem (\ref{50}) is non-degenerate.
\md 

Our second  example is   the problem 
\beq
\lbl{53}
 \p(u'(r))'+\frac{n-1}{r} \p(u'(r))  +b(r)
u^q=0, \s \s r \in (0,1),  \s 
 u =0 \s \mb{for $r=1$}. 
\eeq

By a similar proof, we establish the following result. 
\begin{thm}\lbl{thm:12}
Let $q$ satisfy  (\ref{49}), and assume that $b(r)$ satisfies
\[
 b(r)>0, \s  \s b'(r)<0 \s \s \mb{for } \, r \in (0,1) \,,
\]
\[
\mbox{the functions $\frac{r b'(r)}{b(r)}$ and $r b'(r)$ are non-increasing on $(0,1)$} \,.
\]
Then any positive solution of the problem (\ref{50}) is non-degenerate.
\end{thm}

We consider next an autonomous problem ($f=f(u)$)
\beq
\lbl{54}
\p(u'(r))'+\frac{n-1}{r} \p(u'(r))  +\la f(u)=0, \s \s r \in (0,1),  \s 
 u =0 \s \mb{for $r=1$}, 
\eeq
depending on a positive parameter $\la$. This problem is scaling invariant, which makes it easy to prove the following lemma, see B. Franchi, E. Lanconelli  and J.  Serrin \cite{F}. (Recall that the condition (\ref{14c}) implies that  $u'(r)<0$, so that $u(0)$ gives the maximum value of solution.)

\begin{lma}\lbl{lma:12}
 Assume that the condition (\ref{14c}) holds. Then the maximum value of the positive solution of (\ref{54}),  $u(0)=\al$, uniquely identifies the solution pair $(\la,u(r))$ (i.e., there is at most one $\la$, with at most one solution $u(r)$, so that $ u(0)=\al$). 
\end{lma}

We now turn to the a priori estimates. Several people have already observed that the  recent Liouville type result of   J. Serrin and H. Zou \cite{SZ}, implies a priori estimates of B. Gidas and J. Spruck type \cite{GS}, provided that the point of maximum of solution is bounded away from the boundary, see \cite{AP1}, \cite{FR}. In our case, positive solutions take their maximum at the origin, so that we have the following lemma, whose standard proof we sketch for completeness. 
\begin{lma}\lbl{lma:14}
For the problem (\ref{i1}) assume that the conditions of Lemma \ref{lma:01} hold, and 
\[
\lim _{u \ra \infty} \frac{f(r,u)}{u^q} =b(r) \,,
\]
with $\min (1,p-1)<q<\frac{np-n+p}{n-p}$, and $b(0)>0$. Then there exists a constant $c$, such any positive solution of the problem (\ref{i1}) satisfies
\[
|u|_{L^{\infty}}<c \,.
\]
\end{lma}

\pf
As we have already mentioned, the maximum of any positive solution occurs at $r=0$, by the Lemma \ref{lma:01}. If there is a sequence of unbounded solutions, their rescaling tends to a non-trivial solution of 
\[
\Delta _p u+b(0)u^q=0 \s \mbox{in $R^n$}, 
\]
which is impossible by the result of   J. Serrin and H. Zou \cite{SZ}, see \cite{AP1} and \cite{FR} for more details.
\epf

\md 

We can now prove existence and uniqueness results for model equations.

\begin{thm}\lbl{thm:16} 
Assuming the conditions of Theorem \ref{thm:12}, the problem (\ref{53}) has a unique positive solution.
\end{thm}

\pf
Assume first that $b(r) \equiv 1$. It is known that there exists a  positive solution (proved by using  variational methods). 
To see that the solution is unique, we consider a family of problems 
\beq
\lbl{56}
 \p(u'(r))'+\frac{n-1}{r} \p(u'(r))  +\la
u^q=0, \s \s r \in (0,1),  \s 
 u =0 \s \mb{for $r=1$},  
\eeq
depending on a positive parameter $\la$. Since by Theorem \ref{thm:12}, positive solutions of (\ref{56}) are non-degenerate, we can continue our solution at $\la =1$, for both increasing and decreasing $\la$ on a solution curve, in the Banach space $X$, defined on p. 382 of A. Aftalion and F. Pacella \cite{AP1},   and this solution curve does not admit any turns. (The space $X$, which is a subspace of $C^1(B)$, was originally defined in F. Pacard and T. Rivi\`{e}re \cite{PR}. A. Aftalion and F. Pacella \cite{AP1} showed that one may apply the implicit function theorem, when working in $X$.) By Lemma \ref{lma:00}, solutions stay positive for all $\la$. By scaling, we see that the maximum value $u(0,\la ) \ra 0$, as $\la \ra \infty$, and $u(0,\la ) \ra \infty$, as $\la \ra 0$ (see Lemma \ref{lma:14}). Since this solution curve covers all possible values of $u(0,\la )$, it follows by Lemma \ref{lma:12} that there is only one solution curve. We conclude uniqueness of positive solution at all $\la$, in particular at $\la =1$.

\md 
Turning to the general case, we consider a family of problems
\beq
\lbl{57}
 \p(u'(r))'+\frac{n-1}{r} \p(u'(r))  +
b^{\theta}(r) u^q=0, \s \s r \in (0,1),  \s 
 u =0 \s \mb{for $r=1$},  
\eeq
depending on a  parameter $0 \leq \theta \leq 1$. It is easy to check that the conditions of Theorem \ref{thm:12} hold for all $\theta$.  Using the Theorem \ref{thm:12} and  Lemma \ref{lma:14}, we can continue the solutions between $\theta=0$ and $\theta=1$, to prove existence of positive solutions to 
the problem (\ref{53}). Continuing backwards, between $\theta=1$ and $\theta=0$, we conclude the uniqueness of positive solution to 
the problem (\ref{53}).
\epf

\begin{thm}\lbl{thm:17} 
Assuming the conditions of Theorem \ref{thm:10}, the problem (\ref{50}) has a unique positive solution.
\end{thm}

\pf
We consider a family of problems
\[
\p(u'(r))'+\frac{n-1}{r} \p(u'(r)) -\theta a(r) u^{p-1}+
b(r) u^q=0, \s \s r \in (0,1),  \s 
 u =0 \s \mb{for $r=1$},
\]
depending on a  parameter $0 \leq \theta \leq 1$. When $\theta=0$, we have a unique positive solution, by the preceding result. Arguing as before, we conclude  existence of positive solutions at $\theta =1$.
\epf

\section{A more general result in the one-dimensional case}
\setcounter{equation}{0}
\setcounter{thm}{0}
\setcounter{lma}{0}
We give   an exact description of the curve of positive solutions  for the $p$-Laplace problem in the one-dimensional case 
\beq
\lbl{1}
\varphi(u'(x))'+\la f(u(x))=0 \s\s \mb{for $-1<x<1$}, \s u(-1)=u(1)=0 \,,
\eeq
where $\p (t)=t|t|^{p-2}$, $p>1$, $\la$ a positive parameter, and the function $f(u) \in C^1(\bar R_+)$ satisfies
\beq
\lbl{2}
f(u)<0, \s \mbox{for $u \in (0,\gamma)$}, \s f(u)>0, \s \mbox{for $u >\gamma$} \,,
\eeq
for some $\gamma >0$, and 
\beq
\lbl{3}
f'(u) -(p-1) \frac{f(u)}{u} >0, \s \mbox{for $u >\gamma$} \,.
\eeq
This problem is autonomous, so that it can be posed on any interval. Here we chose the interval $(-1,1)$ for convenience, related to the symmetry of solutions.
\md 

The following lemma gives a known condition for the Hopf's lemma to hold, see J.L. V\'{a}zquez \cite{V}. Moreover, in the present ODE case a completely elementary proof, using the Gronwall's lemma, can be easily given.

\begin{lma}\lbl{lma:0}
Assume that $f(u) \in C^1(\bar R_+)$ satisfies $f(0)=0$, and
\beq
\lbl{14}
-f(s)<c_0 s^{p-1}, \s \mbox{for small $s>0$, and some $c_0>0$} \,.
\eeq
If $u(x)$ is a positive solution of (\ref{1}), then
\[
u'(1)<0 \,.
\]
\end{lma}

We shall need the linearized problem for (\ref{1})
\beq
\lbl{4}
\s\s \s \left( \varphi'(u'(x))w'(x) \right)'+\la f'(u(x))w(x)=0 \s\s \mb{for $-1<x<1$}, \s w(-1)=w(1)=0 \,.
\eeq

We call a solution $u(x)$ of (\ref{1}) {\em non-degenerate}, if the corresponding problem (\ref{4}) admits only the trivial solution $w=0$, otherwise $u(x)$ is a {\em degenerate} solution. We have the following precise description of the solution curve. The uniqueness part was proved previously by J. Cheng \cite{C}.

\begin{thm}\lbl{thm:1}
Assume that  the conditions (\ref{2})  and (\ref{3}) hold, and $\lim _{u \ra \infty} \frac{f(u)}{u^q}>0$ for some $q>\min (1,p-1)$. In case $f(0)=0$, we assume additionally that (\ref{14}) holds. Then the problem (\ref{1}) has at most one positive solution for any $\la $, while  for $\la $ sufficiently small, the problem (\ref{1}) has a unique positive solution. Moreover,  all positive solutions of (\ref{1}) are non-degenerate, they lie on a unique solution curve, extending for $0<\la  \leq \la _0 \leq \infty$. 
\md 

In case $f(0)<0$, we have $\la _0 <\infty$, there are no turns on the solution curve, with $u(0,\la ) \ra \infty$ as $\la \ra 0$, and $u'(\pm 1,\la _0)=0$. I.e., for $\la >\la _0$ there are no positive solutions, while for $0<\la  \leq \la _0$, there is a unique positive non-degenerate solution.
\md 

 In case $f(0)=0$, we have  $\la _0 =\infty$, there are no turns on the solution curve, with $u(0,\la ) \ra \infty$ as $\la \ra 0$, and $\lim _{\la \ra \infty} u(x,\la)=0$, for $x \ne 0$. I.e., for all $\la >0$, there is a unique positive non-singular solution.
\end{thm}

What is remarkable here is that essentially no restrictions are placed on $f(u)$ in the region where it is negative. We shall present the proof, after a series of lemmas, most of which hold under considerably more general conditions.

\begin{lma}\lbl{lma:1}
Assume that $f(u) \in C^1(\bar R_+)$. Then any positive solution of (\ref{1}) is an even function, with $u'(x)<0$ for $x > 0$.
\end{lma}

\pf
If $\xi>0$ is a critical point of $u(x)$, i.e., $u'(\xi)=0$, then the solutions $u(x)$ and $u(2 \xi -x)$ have the same Cauchy data, and so by uniqueness for initial value problems, the graph of $u(x)$ is symmetric with respect to $\xi$, which is impossible. It follows that $x=0$ is the only critical point, the point of global maximum.
\epf

\begin{lma}\lbl{lma:2}
Assume that $f(u) \in C^1(\bar R_+)$, while $u(x)$ and $w(x)$ are any solutions of (\ref{1}) and (\ref{4}) respectively. Then
\beq
\lbl{5}
\varphi'(u') \left(u''w-u'w' \right) =constant, \s \mbox{for all $x \in [-1,1]$} \,.
\eeq
\end{lma}

\pf
Just differentiate the left hand side of (\ref{5}), and use the  equations (\ref{1}) and (\ref{4}).
\epf

\begin{lma}\lbl{lma:3}
Assume that $f(u) \in C^1(\bar R_+)$. If the linearized   problem (\ref{4}) admits a non-trivial solution, then it does not change sign on $(-1,1)$, i.e., we may assume that $w(x)>0$ on $(-1,1)$.
\end{lma}

\pf
Assume on the contrary that $w(x)$ vanishes on $(0,1)$ (the case when $w(x)$ vanishes on $(-1,0)$ is similar). Let $\xi >0$ be the largest root of $w(x)$ on $(0,1)$. We may assume that $w(x)>0$ on $(\xi,1)$, and then $w'(\xi)>0$. Evaluating the expression in (\ref{5}) at $x=\xi$, and at $x=1$,
\[
\varphi'(u'(\xi)) u'(\xi)w'(\xi)=\varphi'(u'(1)) u'(1)w'(1) \,.
\]
The quantity on the left is negative, while the one on the right is non-negative, a contradiction.
\epf 

Given $u(x)$ and $w(x)$, solutions of (\ref{1}) and (\ref{4}) respectively, we again consider the following function, motivated by M. Tang \cite{T}
\[
T(x)=x \left[ (p-1) \varphi(u'(x))w'(x)+\la f(u(x))w(x) \right]-(p-1) \varphi(u'(x))w(x) \,.
\]

The following lemma is proved by a direct computation.

\begin{lma}\lbl{lma:4}
Assume that $f(u) \in C^1(\bar R_+)$, while $u(x)$ and $w(x)$ are any solutions of (\ref{1}) and (\ref{4}) respectively. Then
\beq
\lbl{6}
T'(x)=p \la f(u(x))w(x) \,.
\eeq
\end{lma}

By Lemma \ref{lma:1}, $u(0)$ gives the maximum value of the solution. If we now assume that $f(u)$ satisfies the condition (\ref{2}), then it is easy to see that $f(u(0))>0$ (otherwise we get a contradiction, multiplying (\ref{1}) by $u$, and integrating over $(-1,1)$). So that we can find a point $x_0 \in (0,1)$, such that
\[
u(x_0)=\gamma \,.
\]
(I.e., $f(u(x))>0$ on $(0,x_0)$, and $f(u(x))<0$ on $(x_0,1)$.) Define
\[
q(x)=(p-1)(1-x) \varphi(u'(x))+\varphi'(u'(x)) u(x) \,.
\]

\begin{lma}\lbl{lma:5}
Assume that $f(u) \in C^1(\bar R_+)$ satisfies the condition (\ref{2}). Then
\[
q(x_0)<0 \,.
\]
\end{lma}

\pf
We have $(p-1) \varphi(t)=t\varphi'(t)$, and so we can rewrite
\[
q(x)=\varphi'(u'(x)) \left[(1-x)u'(x)+u(x) \right] \,.
\] 
Since $\varphi'(t)>0$ for all $t \ne 0$, it suffices to show that $r(x) \equiv (1-x)u'(x)+u(x)<0$ on $[x_0,1)$. 
We have $r(1)=0$, and 
\[
r'(x)=(1-x)u''(x)=-(1-x) \la \frac{f(u(x))}{\varphi'(u'(x))}>0, \s \mbox{on $(x_0,1)$} \,,
\]
and the proof follows.
\epf

\begin{lma}\lbl{lma:6}
Assume that $f(u) \in C^1(\bar R_+)$ satisfies the conditions (\ref{2}) and (\ref{3}), while $u(x)$ and $w(x)$ are any solutions of (\ref{1}) and (\ref{4}) respectively. Then
\beq
\lbl{7}
(p-1)w(x_0)\varphi(u'(x_0))-u(x_0)w'(x_0) \varphi'(u'(x_0))>0 \,,
\eeq
which implies, in particular,
\beq
\lbl{8}
w'(x_0)<0 \,.
\eeq
\end{lma}

\pf
By a direct computation
\[
\left[(p-1)w(x)\varphi(u'(x))-u(x)w'(x) \varphi'(u'(x)) \right]'=\la \left[ f'(u) -(p-1) \frac{f(u)}{u} \right]uw\,.
\]
The quantity on the right is positive on $(0,x_0)$, in view of our conditions, and Lemma \ref{lma:3}. Integrating over $(0,x_0)$, we conclude (\ref{7}).
\epf

We have all the pieces in place for the following crucial lemma. 

\begin{lma}\lbl{lma:7}
Under the conditions (\ref{2}) and (\ref{3}), any positive solution of (\ref{1}) is non-degenerate, i.e., the corresponding linearized problem (\ref{4}) admits only the trivial solution.
\end{lma}

\pf
Assuming the contrary, let $w(x)>0$ on $(-1,1)$ be  a solution of (\ref{4}) (see Lemma \ref{lma:3}). 
By Lemma \ref{lma:2}, and since $f(u(x_0))=0$,
\beq
\lbl{10}
\varphi'(u'(1)) u'(1)w'(1)=\varphi'(u'(x_0)) u'(x_0)w'(x_0)=(p-1) \varphi(u'(x_0)) w'(x_0) \,.
\eeq
Integrating (\ref{6}) over $(x_0,1)$, we have $T(1)-T(x_0)=
p \la \int_{x_0}^1 f(u(x))w(x) \, dx<0$, i.e.,
\[
L \equiv (p-1)\varphi(u'(1))w'(1)- (p-1) x_0 \varphi(u'(x_0))w'(x_0)+(p-1) \varphi(u'(x_0))w(x_0)<0 \,.
\]
On the other hand, using (\ref{10}), then (\ref{7}), followed by (\ref{8}) and Lemma \ref{lma:5}, we have
\beqa \nonumber
& L=(p-1)\varphi(u'(x_0))w'(x_0)-(p-1)x_0 \varphi(u'(x_0))w'(x_0)+(p-1)\varphi(u'(x_0))w(x_0) \\ \nonumber
&  >(p-1)\varphi(u'(x_0))w'(x_0)-(p-1)x_0 \varphi(u'(x_0))w'(x_0)+u(x_0)w'(x_0) \varphi'(u'(x_0)) \\\nonumber
& =w'(x_0)q(x_0)>0 \,,\nonumber
\eeqa
a contradiction.
\epf

This lemma implies that there are no turns on the solution curves. For a more detailed description of the bifurcation diagram, we need to establish some additional properties of solutions. Denoting $F(u)=\int _0^u f(t) \, dt$, one shows by differentiation that
\beq
\lbl{12}
E(x) \equiv \frac{p-1}{p} |u'(x)|^p+\la F(u(x))=constant, \s \mbox{for all $x \in [-1,1]$}.
\eeq
It is customary to think of $E(x)$ as ``energy".
\md

We shall also need the  following well-known lemma, see e.g. P. Korman \cite{K}.

\begin{lma}\lbl{lma:10}
 The maximum value of the positive solution of (\ref{1}),  $u(0)=\al$, uniquely identifies the solution pair $(\la,u(x))$ (i.e., there is at most one $\la$, with at most one solution $u(x)$, so that $ u(0)=\al$). 
\end{lma}

\noindent
{\bf Proof of the Theorem \ref{thm:1}} We begin by showing that the problem (\ref{1}) has solutions for some $\la >0$. Observe that our condition (\ref{3}) implies that the function $\frac{f(u)}{u^{p-1}}$ is increasing, i.e., $f(u) \ra \infty$ for large $u$. Let $\theta>\gamma$ be the point where $\int_0^{\theta} f(t) \, dt=0$. Consider the solutions of  the initial value problem
\[
\varphi(u')'+ f(u)=0, \s u(0)=\al, \s u'(0)=0 \,.
\]
If $\al >\theta$, then the conservation of energy formula (\ref{12}) implies that $\frac{p-1}{p} |u'(x)|^p+ F(u(x))=F(\al )>0$, and so this solution must  vanish (for the first time) at some $r>0$. Rescaling, we get a positive solution of (\ref{1}) at $\la =r^2$.
\md 

We now use the implicit function theorem to continue solutions for both decreasing and increasing $\la$ (working, as before, in the space $X$ from \cite{AP1}, or in the framework of B.P. Rynne \cite{R}), and in both directions there are no turns by Lemma \ref{lma:7}. For decreasing $\la$, we have $u(0,\la) \ra \infty$, since otherwise the solution curve would have no place to go (it cannot go to zero, because $f(u)$ is negative near zero). By Lemma \ref{lma:14}, $u(0,\la) \ra \infty$ as $\la \ra 0$. For increasing $\la$, $u(0,\la)$ is decreasing, in view of Lemma \ref{lma:10}, and hence the oscillation of the solutions is decreasing. Consider first the case $f(0)=0$. Then positive solutions continue for all $\la$, since $u_x(1,\la)<0$, by Lemma \ref{lma:0}.
In case $f(0)<0$,  the opposite is true, i.e.,  we have $u_x(1,\la _0)=0$ at some $\la _0$, and the solution becomes sign-changing for $\la > \la _0$, since otherwise we shall get arbitrarily large oscillations of $u(x, \la)$ for large $\la$, by integrating the equation (\ref{1}) over any interval $(x,1)$.
\epf

\end{document}